\RequirePackage{ifpdf}
\ifpdf 
\documentclass[pdftex]{sigma}
\else
\documentclass{sigma}
\fi

\usepackage{amscd}

\usepackage{amstext,eucal}
\usepackage{bbold}

\renewcommand{\phi}{\varphi}

\newcommand{\cA}{{\mathcal{A}}}

\newcommand{\cW}{{\mathcal{W}}}

\newcommand{\cE}{{\mathcal{E}}}

\newcommand{\cK}{{\mathcal{K}}}

\newcommand{\cT}{{\mathcal{T}}}

\newcommand{\cS}{{\mathcal{S}}}
\newcommand{\cD}{{\mathcal{D}}}
\newcommand{\cP}{{\mathcal{P}}}

\newcommand{\Rho}{{\mbox{\sf P}}}

\newcommand{\RR}{\mathbb{R}}

\newcommand{\CC}{\mathbb{C}}

\newcommand{\SO}{\mathrm{SO}}
\newcommand{\Spin}{\mathrm{Spin}}

\begin{document}

\allowdisplaybreaks

\renewcommand{\thefootnote}{$\star$}

\renewcommand{\PaperNumber}{079}

\FirstPageHeading

\ShortArticleName{About Twistor Spinors with Zero in Lorentzian Geometry}

\ArticleName{About Twistor Spinors with Zero\\ in Lorentzian Geometry\footnote{This paper is a
contribution to the Special Issue ``\'Elie Cartan and Dif\/ferential Geometry''. The
full collection is available at
\href{http://www.emis.de/journals/SIGMA/Cartan.html}{http://www.emis.de/journals/SIGMA/Cartan.html}}}

\Author{Felipe LEITNER}

\AuthorNameForHeading{F. Leitner}

\Address{Universit{\"a}t Stuttgart, Institut f{\"u}r Geometrie und Topologie,
Fachbereich  Mathematik,\\
 Pfaffenwaldring 57,
D-70550 Stuttgart, Germany}
\Email{\href{mailto:leitner@mathematik.uni-stuttgart.de}{leitner@mathematik.uni-stuttgart.de}}
\URLaddress{\url{http://www.igt.uni-stuttgart.de/LstDiffgeo/Leitner/}}

\ArticleDates{Received April 06, 2009, in f\/inal form July 10, 2009;  Published online July 28, 2009}

\Abstract{We describe the local conformal geometry of a Lorentzian spin manifold $(M,g)$ admitting a twistor spinor $\phi$ with zero.
Moreover, we describe the shape of the zero set of $\phi$. If $\phi$ has isolated zeros then the metric $g$ is locally conformally
equivalent to a static monopole. In the other case the zero set consists of null geodesic(s) and  $g$ is locally conformally
equivalent to a Brinkmann metric. Our arguments utilise tractor calculus in an essential way. The Dirac current of $\phi$, which is
a conformal Killing vector f\/ield, plays an important role for our discussion as well.}


\Classification{53C27; 53B30}

\section{Introduction}\label{Sec1}

There are two conformally covariant f\/irst order partial dif\/ferential equations for spinors, the {\it Dirac equation}
and the {\it twistor equation}. The twistor equation is an  overdetermined  system of PDE's, which was f\/irst introduced by R.~Penrose in
the context of General Relativity (cf.~\cite{PR86}). We are interested in this article in solutions of the twistor equation with zeros on Lorentzian spin manifolds.

Let $(M^n,g)$ be a time-oriented  Lorentzian manifold of dimension $n\geq 3$ with spin structure and spinor bundle $\cS$. The twistor equation for spinors
$\phi\in\Gamma(\cS)$ is given by
\[
\nabla^\mathcal{S}_X\phi +\frac{1}{n}X\cdot \cD\phi=0\qquad\mbox{for\ any}\quad X\in TM
\]
(cf.~Section~\ref{Sec2}). Its solutions are called {\it twistor spinors} or {\it conformal Killing spinors}. We prefer the latter notion in this article
(apart from the title and abstract).

There are many geometric structure results for conformal Killing spinors in the literature (cf.\ e.g.~\cite{Lic88,Lew91,BFGK91,BL04} and the surveys~\cite{B00,B08}).
Especially, the
case
of
Riemannian
geometry
is
well
studied.
In particular, it is known that a Riemannian metric admitting a  conformal Killing spinor with zero is conformally equivalent to an
asymptotically locally Euclidean and Ricci-f\/lat  metric. This result is achieved by conformally rescaling the metric with the norm of the conformal
Killing spinor.
The aim of this article is to derive an analogous result in Lorentzian geometry. In particular, we study
in the following the local conformal geometry and the shape of the
zero set of conformal Killing spinors on Lorentzian manifolds. However, an analogous treatment of the problem as in Riemannian geometry is not possible for various
reasons.
Instead we employ in this article methods of tractor
calculus and its spinorial counterpart -- twistor calculus (cf.~\cite{Tho34,BEG94,L07b}). This approach and its results demonstrate the force of tractor (and twistor)
calculus
in a particular instance.

Solutions of the twistor equation with zeros on Lorentzian manifolds are interesting, in particular, for the following reason. A conformal Killing spinor
$\phi$ is always accompanied by a~conformal Killing vector f\/ield $V_\phi$, the so-called {\it Dirac current}. If the conformal Killing spinor~$\phi$ admits a~zero, then
the  Dirac current $V_\phi$ is an  essential conformal  Killing vector f\/ield. In particular, the local f\/low of $V_\phi$ consists of conformal
transformations, which do not preserve any metric in the underlying conformal class! A conjecture by A. Lichnerowicz states that any compact Lorentzian
manifold with an essential conformal transformation group is conformally f\/lat (cf.~\cite{DAG91}). In the analytic category Ch.~Frances proved in \cite{F07}
that any Lorentzian manifold, admitting a causal conformal Killing vector with a zero, is conformally f\/lat. On the other hand, we have constructed
in \cite{L07a} a Lorentzian spin manifold admitting a conformal Killing spinor  with isolated  zero, which is not conformally f\/lat in the spacelike neighbourhood
of the isolated zero. This is possible, since the construction in \cite{L07a} is not analytic, but only $C^1$-dif\/ferentiable.

In Section \ref{Sec2} we recall the twistor equation in the classical way as partial dif\/ferential equation for spinor f\/ields. In Section \ref{Sec3}
we introduce conformal tractor and twistor calculus via Cartan geo\-metry.   This approach allows us to reformulate the Penrose
twistor equation in terms of parallel twistor f\/ields. Spinors and twistors are accompanied by so-called spinor and twistor squares, respectively.
In Lorentzian geometry the spinorial square of f\/irst degree plays a particular role. This is the so-called Dirac current. For conformal
Killing spinors $\phi$ the Dirac current $V_\phi$ is a~conformal Killing vector f\/ield. On the other hand, to $\phi$ and $V_\phi$
there corresponds a parallel twistor with adjoint tractor. All this is explained in Section~\ref{Sec4}. The tractor/twistor picture is
the essential tool for the geometric
description
of the underlying Lorentzian conformal structure. Brinkmann spaces and static monopole solutions are the possible underlying geometries. This is stated
in Proposition~\ref{prop1}
of
Section
\ref{Sec5}. Finally, in Section \ref{Sec6} we discuss the shape of the zero set ${\rm zero}(\phi)$ and the singular set ${\rm sing}(\phi)$
of a~conformal Killing spinor $\phi$. Either isolated zeros or null geodesics are possible. These results about the shape are based on investigations
of  \cite{L99} (see also \cite{L01}).

\section{The twistor equation for spinors}\label{Sec2}

We brief\/ly recall here the twistor equation for spinors (in Lorentzian geometry),
which is most notably a conformally covariant, overdetermined partial dif\/ferential equation of f\/irst order on manifolds (cf.~\cite{PR86,BFGK91}).

Let  $\SO(r,s)$ denote the  special orthogonal group, which acts by the standard representation on  the Euclidean space
$\RR^{r,s}=(\RR^n,(\cdot,\cdot)_{r,s})$ of dimension $n=r+s$ with (indef\/inite) scalar product of signature $(r,s)$ (where $r$ denotes the number of
timelike vectors in an orthonormal basis of~$\RR^{r,s}$). The  spin group of signature
$(r,s)$
is denoted by  $\Spin(r,s)$ and covers  $\SO(r,s)$ twice by a~group homomorphism
\[\lambda: \  \Spin(r,s)\to \SO(r,s) .\]
Then let $\mathbb{\Delta}_{r,s}$ denote the  standard complex spinor module with spin representation $\rho=\rho_{r,s}$ (cf.\ e.g.~\cite{LM:89,Bau81,BFGK91}).
We denote the connected components of the identity elements in  $\SO(r,s)$ and  $\Spin(r,s)$ by  $\SO_o(r,s)$ and  $\Spin_o(r,s)$,
respectively.

Now let $(M^n,g)$ be a connected and time-oriented Lorentzian spin manifold (i.e.~$r=1$) of dimension $n\geq 3$.
The principal $\SO_o(1,n-1)$-bundle of time- and space-oriented orthonormal frames on $M$ is denoted by $SO(M)$. We f\/ix a spin structure
$(\Spin(M),\pi)$ on $M$, that is a~$\lambda$-reduction of  $SO(M)$ to the spin group $\Spin_o(1,n-1)$. The (complex) spinor bundle is given by
\[\mathcal{S}=\Spin(M)\times_{\rho}\mathbb{\Delta}_{1,n-1} . \]
We denote the set of smooth sections in $\mathcal{S}$ by $\Gamma(\mathcal{S})$  and refer to them as {\it spinor f\/ields}
or just {\it spinors} on $(M,g)$.
The  spinor bundle $\mathcal{S}$ is equipped with an ($\Spin_o(1,n-1)$-invariant)
indef\/inite Hermitian product $\langle\cdot,\cdot\rangle_{\cS}$.
Moreover, the Levi-Civita connection $\omega^g$ on  $SO(M)$ lifts via $\pi$ to a~connection on the spin frame bundle
$\Spin(M)$, which in turn induces a covariant derivative on~$\mathcal{S}$, the so-called {\it spinor derivative}
\[\nabla^{\mathcal{S}}: \ \Gamma(\mathcal{S})\to \Omega^1(M)\otimes \Gamma(\mathcal{S}) .\]
The spinor derivative $\nabla^{\mathcal{S}}$ preserves the  Hermitian  product  $\langle\cdot,\cdot\rangle_{\cS}$.

The bundle $T^*M\otimes\cS$ of spinor-valued $1$-forms decomposes into the direct sum $\mathcal{S}\oplus\mathcal{K}$,
where $\mathcal{K}$ is the kernel of the Clif\/ford multiplication of $T^*M$ on $\mathcal{S}$. The superposition of
$\nabla^{\mathcal{S}}$ with the projection $\pi_\cS$ onto  $\mathcal{S}$ yields the {\it Dirac operator} $\mathcal{D}^g:\Gamma(\mathcal{S})\to\Gamma(\mathcal{S})$
on $(M,g)$.
The projection $\mathcal{P}^g:=\pi_\cK\circ\nabla^{\mathcal{S}}:\Gamma(\cS)\to\Gamma(\mathcal{K})$ is the so-called {\it Penrose
operator}. The (non-trivial) elements of the kernel of the Penrose operator $\mathcal{P}$ are called {\it twistor spinors}, or alternatively, {\it conformal Killing
spinors}.
(Apart from the title we use in this article the latter notation.)
It follows easily from the  def\/inition that a spinor f\/ield $\phi\in\Gamma(\cS)$ is a conformal Killing spinor if and only if
$\phi$ satisf\/ies the  {\it Penrose twistor equation}
\begin{equation}\label{eq1}
\nabla^\mathcal{S}_X\phi +\frac{1}{n}X\cdot \cD\phi=0\qquad\mbox{for\ any}\quad X\in TM,
\end{equation}
where $X\cdot \cD\phi$ denotes the Clif\/ford product of $X$ with the spinor $\cD\phi$.

In this  article,  we are interested in Lorentzian spin manifolds  $(M^n,g)$, $n\geq 3$, which admit a~non-trivial solution $\phi$ of (\ref{eq1})
such that the zero set ${\rm zero}(\phi)=\{p\in M|\, \phi(p)=0\}$ is non-empty in $M$. In particular, we are interested in the shape of ${\rm zero}(\phi)$
and the geometric properties of $(M^n,g)$ (locally outside of the zero set ${\rm zero}(\phi)$).

The conformal class of a metric $g$ on $M$ consists of all those metrics $\tilde{g}$, which dif\/fer from $g$ only by multiplication
with a positive function on $M$, i.e., $\tilde{g}=e^{2f}g$ for some $f\in C^\infty(M)$.
It is a well known fact that $\cD$ and $\mathcal{P}$ are conformally covariant dif\/ferential operators.
In case of the Penrose operator $\mathcal{P}$ this means \begin{equation}\label{eq1.5}\cP^{\tilde{g}}=e^{-f/2}\circ \cP^g\circ e^{-f/2}\end{equation}
for $\tilde{g}=e^{2f}g$, i.e., the bidegree of  $\mathcal{P}$ is $(-1/2,1/2)$.   In particular, the kernel
of  $\mathcal{P}$ is a conformal invariant of $(M,c)$. The framework of conformal Cartan geometry provides a tool, in particular, for
the  invariant construction of the  Penrose twistor equation. We aim to explain this next. This approach will be
the key to our discussion of conformal Killing spinors with zeros on  Lorentzian manifolds.

\section{Tractors and twistors}\label{Sec3}

We introduce here the {\it standard tractor bundle} $\cT$ with tractor connection $\nabla^\cT$ of conformal (Lorentzian) geometry via Cartan geometry (cf.~\cite{Tho34,BEG94}).
The spinorial counterpart is the so-called {\it twistor bundle} $\cW$ with connection (cf.\ e.g.~\cite{L07b}).
We then argue that conformal Killing spinors correspond uniquely to parallel sections of $\cW$ (cf.~\cite{PM72,Fri76,BFGK91}).

The projective null cone $\mathbb{P}L$  of the Euclidean space $\RR^{2,n}$ of signature $(2,n)$
serves as the standard model of conformal Lorentzian geometry. As a homogeneous space,  $\mathbb{P}L$ is given by $\SO_o(2,n)/P$,
where $P\subset \SO_o(2,n)$ is the parabolic subgroup, which stabilises a null line in $\RR^{2,n}$. The group $\SO_o(2,n)$ acts naturally by conformal
transformations on $\mathbb{P}L$.

Now let $(M^n,g)$ be a time- and space-oriented Lorentzian manifold of dimension $n\geq 3$. We denote the conformal class
of $g$ on $M$ by $c=[g]$. Thus we obtain a {\it conformal Lorentzian manifold} $(M^n,c)$.
The  conformal structure $c=[g]$ on $M$ is equivalently described by a  principal $P$-bundle $P(M)$ over $M$
with canonical  Cartan connection
\[
\omega^c: \ P(M)\to \frak{so}(2,n)  ,
\]
where $\frak{so}(2,n)$ denotes the orthogonal Lie algebra of signature $(2,n)$.
The canonical Cartan connection $\omega^c$ is uniquely determined by the condition $\partial^*\kappa=0$, where $\kappa:P(M)\to {\rm Hom}(\Lambda^2\RR^n,\frak{so}(2,n))$
denotes the curvature function of  $\omega^c$ and $\partial^*$ is the {\it Kostant-codifferential} (cf.~\cite{Kob72,CSS97a}).

The associated vector bundle
\[\cT\ :=\ P(M)\times_\iota \RR^{2,n}, \]
where $\iota$ denotes the standard representation of $P\subset\SO_o(2,n)$ on  $\RR^{2,n}$, is called the {\it standard tractor bundle} of $(M,c)$ (cf.~\cite{Tho34,BEG94}).
The tractor bundle $\cT$ is equipped with an invariant metric $\langle\cdot,\cdot\rangle_{\cT}$, and the Cartan connection $\omega^c$ induces a
$\langle\cdot,\cdot\rangle_{\cT}$-compatible covariant derivative
$\nabla^\cT$ on $\cT$. Moreover, $\cT$  admits a $P$-invariant f\/iltration
\[\cT\ \supset\ \cT^0\ \supset\ \cT^1,\]
where $\cT^1$ is the subbundle of $\cT$, which corresponds to the $P$-stable null line of the standard representation on $\RR^{2,n}$, and $\cT^0$
is the $\langle\cdot,\cdot\rangle_{\cT}$-orthogonal bundle to $\cT^1$  in $\cT$.

Note that any metric $g\in c$ in the conformal class on $M$ serves as a so-called {\it Weyl structure} in the sense of \cite{CS03} and reduces the structure
group
of
$P(M)$ to $\SO_o(1,n-1)$. This reduction of the structure group with respect to $g\in c$ causes a splitting of the f\/iltration of $\cT$
into a direct sum of the form $\RR\oplus TM\oplus\RR$, where the f\/irst summand corresponds to $\cT/\cT^0$ and the last summand
to  $\cT^1$. We set $s_-:=(1,0,0)$ and $s_+:=(0,0,1)$. Any section in  $\cT\to M$  is given with respect to $g\in c$  by a triple $(a,V,b)=as_-+V+bs_+$
with $a,b\in C^\infty(M)$ and $V\in\frak{X}(M)$.
The tractor metric satisf\/ies $\langle(a,V,b),(\tilde{a},\tilde{V},\tilde{b})\rangle_\cT=a\tilde{b}+b\tilde{a}+g(V,\tilde{V})$.

In case $M$ is a spin manifold for some (hence all) $g\in c$, the conformal Cartan geometry  $(P(M),\omega^c)$ over $M$ lifts to
a conformal spin Cartan geometry  $(\tilde{P}(M),\omega^c)$, where  $\tilde{P}$ in $\Spin_o(2,n)$ is the preimage of $P$ under
$\lambda: \Spin_o(2,n)\to \SO_o(2,n)$, and $\tilde{P}(M)$ is a principal  $\tilde{P}$-bundle.
We call the $\rho_{2,n}$-associated vector bundle
\[
\mathcal{W}=\tilde{P}(M)\times_{\rho}\mathbb{\Delta}_{2,n}
\]
the {\it twistor bundle} of the conformal Lorentzian spin manifold $(M,[g])$ (with f\/ixed spin structure).
The twistor bundle $\mathcal{W}$ is equipped with an indef\/inite Hermitian inner product $\langle\cdot,\cdot\rangle_{\cW}$ and an $\omega^c$-induced
covariant derivative $\nabla^\cW$ preserving $\langle\cdot,\cdot\rangle_{\cW}$. Moreover, we have a f\/iberwise Clif\/ford multiplication $\cdot:\cT\otimes\cW\to\cW$
of tractors with twistors. There is also a natural subbundle~$\cW^{1/2}$ in~$\cW$, that is the bundle of twistors, which are annihilated by Clif\/ford multiplication
with the elements of the real line bundle $\cT^1$. This def\/ines a natural f\/iltration $\cW\supset \cW^{1/2}$ on the twistor bundle, and
the quotient  bundle  $\mathcal{W}/\cW^{1/2}$
is naturally isomorphic to the spinor bundle $\cS[-1/2]:=\cS\otimes\mathcal{E}[-1/2]$ of conformal weight $-1/2$. Here $\mathcal{E}[w]$, $w\in\RR$,
denotes the density bundle, which is the trivial real line bundle over $M$ associated to the representation $|\det|^{w/n}$ of the general linear group
$\mathrm{GL}(n)$.

Any metric $g\in c$ in the conformal class on $M$ reduces the structure
group
of
$\tilde{P}(M)$ to $\Spin_o(1,n-1)$. Accordingly, the f\/iltration $\mathcal{W}\supset\cW^{1/2}$
splits into the direct sum
\begin{equation}\label{eq2}
\mathcal{W}\ \cong_g\ \cS\oplus\cS  ,
\end{equation}
where both, the quotient $\mathcal{W}/\cW^{1/2}$ and the subbundle $\cW^{1/2}$, become isomorphic to the spinor bundle $\cS$ with respect to $g$ on $M$.
Note that the invariant Hermitian product on $\cW$ is given with respect to the $g$-corresponding splitting (\ref{eq2})
by
\[
\langle(\phi,\psi),(\tilde{\phi},\tilde{\psi})\rangle_{\cW}  =  \frac{i}{\sqrt{2}}\big(  \langle \phi, \tilde{\psi}\rangle_{\cS}- \langle
\psi,\tilde{\phi}\rangle_{\cS}\big) .
\]

Now let us consider the most obvious condition $\nabla_\cdot^\cW\Phi=0$ for a twistor $\Phi\in\Gamma(\cW)$. This equation means
that $\Phi$ is a parallel twistor on $(M,c)$ with respect to the canonical connection. It is clear by construction that  this is a  conformally
invariant  condition  for  a twistor $\Phi$.
With respect to any $g\in c$ and
the
induced splitting
$\Phi=(\phi,\psi)$ of (\ref{eq2}) the equation  $\nabla_\cdot^\cW\Phi=0$ is equivalent to the Penrose equation (\ref{eq1}) for $\phi$ and
the equation $\nabla_X^\cS\psi=\frac{1}{\sqrt{2}}\Rho^g(X)\cdot \phi$ for any $X\in TM$,  where $\Rho^g(X)\cdot \phi$ is the Clif\/ford product with the image of
the Schouten operator  $\Rho^g$ on $(M,g)$. The Schouten operator $\Rho^g:TM\to TM$ is def\/ined by
\[\Rho^g=\frac{1}{n-2}\left(\frac{{\rm scal}^g}{2(n-1)}{\rm id}|_{TM}-{\rm Ric}^g\right) ,\]
where ${\rm Ric}^g$ is the Ricci curvature operator and ${\rm scal}^g$ is the scalar curvature of $(M,g)$.
On the other hand, any conformal Killing spinor $\phi\in\Gamma(\cS)$ on $(M,g)$ satisf\/ies the equation
$\nabla_X^\cS\cD\phi=\frac{n}{2}\Rho^g(X)\cdot \phi$.
Thus a parallel twistor $\Phi$ on $(M,c)$ corresponds via the splitting (\ref{eq2})
with respect to some $g\in c$ uniquely to a conformal Killing spinor $\phi\in\Gamma(\cS)$ on $(M,g)$. The correspondence is explicitly given by
\[\phi\quad \stackrel{g}{\longleftrightarrow}\quad \Phi=\bigg( \phi , \frac{\sqrt{2}}{n}\cD^g\phi \bigg)\]
(cf.~\cite{BFGK91,L01}).
Also note that,  with respect to a metric $g\in c$, the Penrose operator $\cP^g$ is the f\/irst dif\/ferential operator in the spinorial BGG sequence,
which maps sections of $\cS[-1/2]=\mathcal{W}/\cW^{1/2}$ to sections of $\cK[1/2]=\cK\otimes \cE[1/2]$
on $(M,c)$ (cf.~\cite{CSS01}). This explains the bidegree in~(\ref{eq1.5}) and, in particular,  conformal Killing spinors naturally have conformal weight~$-1/2$.

\section{The Dirac current and its adjoint tractor}\label{Sec4}

We recall here the notion of the {\it Dirac current} of a spinor $\phi$ on a Lorentzian spin manifold.
On the other hand, any twistor $\Phi$ on a conformal Lorentzian spin
manifold gives rise  in a natural way to an adjoint (resp. $2$-form) tractor.

In general, the spinor module $\mathbb{\Delta}_{r,s}$ can be understood as a square root of the complexif\/ied exterior algebra
$\Lambda_\CC:=\Lambda(\RR^{r,s})\otimes\CC$  of  the
Euclidean space $\RR^{n}$ with scalar product $(\cdot,\cdot)_{r,s}$. To be concrete, as representation spaces of $\Spin(r,s)$ we have
\begin{alignat*}{3}
& \zeta: \ \mathbb{\Delta}_{r,s}\otimes  \mathbb{\Delta}_{r,s} \cong \Lambda_\CC \qquad &&\mbox{for}\ n\ \mbox{even\quad and} & \\
& \zeta: \ \mathbb{\Delta}_{r,s}\otimes  \mathbb{\Delta}_{r,s} \cong \Lambda^{\rm ev}_\CC\cong\Lambda^{\rm odd}_\CC\qquad & &\mbox{for}\ n\
\mbox{odd},&
\end{alignat*}
where $\Lambda^{\rm ev}_\CC$ and  $\Lambda^{\rm odd}_\CC$ denote the exterior products of even and odd degree, respectively.
(The isomorphism  $\Lambda^{\rm ev}_\CC\cong\Lambda^{\rm odd}_\CC$ is realised by the Hodge-$\star$ operator.)
In particular, for any $\phi\in\mathbb{\Delta}_{r,s}$
and $0\leq i\leq n$
we have the  notion of a {\it spinor square} $\zeta^i(\phi)\in \Lambda^{i}_\CC$ of degree $i$, which is by def\/inition the homogeneous component of degree
$i$ of $\zeta(\phi\otimes\phi)$. Recall that the exterior algebras of $\RR^{r,s}$ and its dual space ${\RR^{r,s}}^*$ are canonically identif\/ied
via the metric  $(\cdot,\cdot)_{r,s}$.
Accordingly, we have the dual (complex) $i$-form to  $\zeta^i(\phi)$ on $\RR^{r,s}$, which we denote  by $\zeta_i(\phi)$.

Especially, in the Lorentzian case $r=1$, the (real) vector $\zeta^1(\phi)\in\RR^{1,n-1}$ is explicitly given via the relation
\[
(\zeta^1(\phi),x)_{1,n-1}=-\langle x\cdot \phi,\phi\rangle_{1,n-1}\qquad \mbox{for\ all}  \quad x\in \RR^{1,n-1}  ,
\]
where $\langle \cdot,\cdot\rangle_{1,n-1}$ is the Hermitian product on the spinor module $\mathbb{\Delta}_{1,n-1}$.
The vector $\zeta^1(\phi)$ is causal (i.e.\ null or timelike) and future-directed for any $\phi\neq 0$ (cf.\ e.g.~\cite{Bau81,L01}).
Accordingly, we obtain for any spinor f\/ield $\phi\in\Gamma(\cS)$ on some  time-oriented Lorentzian spin manifold $(M^n,g)$ a vector f\/ield $V_\phi$,
which is uniquely determined by the relation $g(V_\phi,X)=-\langle X\cdot \phi,\phi\rangle_{\cS}$ for all $X\in TM$. We call $V_\phi\in\frak{X}(M)$ the
{\it Dirac current} of the spinor $\phi\in\Gamma(\cS)$ on the Lorentzian manifold $(M,g)$. The Dirac current $V_\phi$ is a causal and future-directed
vector f\/ield, whose zero set coincides with the zero set of $\phi\in\Gamma(\cS)$:
\[
{\rm zero}(V_\phi)= {\rm zero}(\phi)\qquad   \mbox{on}\ \ M .
\]
We denote the corresponding   dual $1$-form to $V_\phi$ via $g$ by $\alpha_\phi\in\Omega^1(M)$.

Now let $(M^n,c)$ be a time-oriented  conformal Lorentzian spin manifold of dimension $n\geq 3$. Recall that we have the tractor bundle $\cT$ with  dual
$\cT^*$
and the
twistor
bundle $\cW$ over
$(M,c)$.
Moreover, we have the adjoint tractor bundle $\cA$ on $(M,c)$, which is the associated vector bundle  to the adjoint representation of the parabolic
subgroup $P$ in
$\mbox{SO}_o(2,n)$:
\[
\cA := P(M)\times_{\rm Ad}  \frak{so}(2,n) .
\]
The adjoint tractor bundle $\cA$ is canonically isomorphic to the bundle $\Lambda^2\cT^*$ of $2$-form tractors on $(M,c)$.
(In the following, we will not distinguish between adjoint and $2$-form tractors.) Note that the squaring $\zeta$
induces for any twistor $\Phi\in\Gamma(\cW)$ a $2$-form tractor $\alpha_\Phi\in\Gamma(\cA)$ on $(M,c)$. This $2$-form  tractor is explicitly given by
the relation
\[
\langle  \alpha_\Phi, A\rangle_\cT  =  -i\langle A\cdot \Phi,\Phi\rangle_\cW\qquad  \mbox{for\ any} \quad A\in\cA,
\]
where  $A\cdot \Phi$ denotes the Clif\/ford product of a $2$-form tractor with a twistor and $\langle\cdot, \cdot\rangle_\cT$ is the tractor metric
extended to~$\cA$.

We have seen in the previous section that with respect to a Lorentzian metric $g\in c$ the tractor bundle $\cT$ and the twistor bundle $\cW$
split into direct sums over $M$. Also the adjoint tractor bundle $\cA$ splits with respect to $g$ into a direct sum, which is
\[ T^*M\oplus \frak{co}(TM)\oplus T^*M,\]
where $\frak{co}(TM)$ denotes the bundle of exterior $2$-forms on $M$ plus multiples of the symmetric $(2,0)$-tensor   $g$.
Equivalently, using the dual tractors $s^\flat_-:=\langle s_-,\cdot\rangle_\cT$ and $s^\flat_+:=\langle s_+,\cdot\rangle_\cT$, any $2$-form tractor $A\in\Gamma(\cA)$ is
given in the form
\[
s_-^\flat\wedge \alpha_-  +  \alpha_0 +  s^\flat_-\wedge s^\flat_+\wedge \alpha_\mp +   s^\flat_+\wedge\alpha_+
\]
with uniquely determined dif\/ferential forms $\alpha_-,\alpha_+\in\Omega^1(M)$, $\alpha_0\in\Omega^2(M)$ and $\alpha_\mp\in C^\infty(M)$.
We call $\alpha_-$ the {\it projecting component} of $A$ (cf.~\cite{LeiNCK}).

Now let $\Phi=(\phi,\psi)$ be some twistor given as a pair of spinors with respect to the splitting~(\ref{eq2}) induced by $g\in c$.
It is natural to ask how the twistor square $\alpha_\Phi$ and the spinor square $\alpha_\phi$ (resp.\ the Dirac current $V_\phi$) of the f\/irst component
of $\Phi$ are related (with respect to $g\in c$). A~straightforward calculation shows $\alpha_\Phi=(\alpha_-,\alpha_0,\alpha_\mp,\alpha_+)$ with
\[
\alpha_-=\alpha_\phi,\qquad \alpha_\mp=\sqrt{2}\cdot \mbox{Re}\langle\phi,\psi\rangle_\cS,\qquad \alpha_+=-\alpha_{\psi},
\]
and the $2$-form $\alpha_0$ on $M$ is determined  by $\alpha_0(X,Y)=\frac{1}{\sqrt{2}} \mbox{Re}\langle X\wedge Y\cdot \phi,\psi\rangle_\cS$ for $X,Y\in TM$.

Now let $\phi\in\Gamma(\cS)$ be a conformal Killing spinor on $(M^n,g)$. We have seen in the previous section that $\phi$ corresponds in a unique way
to a  $\nabla^\cW$-parallel twistor $\Phi\in\Gamma(\cW)$ on the conformal Lorentzian spin manifold $(M,[g])$, which is given with respect to $g$ by
the pair $\Phi=(\phi,\frac{\sqrt{2}}{n}\cD^g\phi)$.
Note that, by construction, the induced $2$-form tractor $\alpha_\Phi$ of  $\Phi$ on   $(M,[g])$ is  parallel with respect to the extended tractor connection
$\nabla^\cT$ on $\cA$.
In \cite{LeiNCK,Lei04a} we have studied the projecting component of $\nabla^\cT$-parallel $q$-form tractors, in general. The discussion shows that the projecting
component
of any  $\nabla^\cT$-parallel $2$-form tractor is a conformal Killing $1$-form $\alpha_-$, i.e., the
symmetric trace-free part of the $(2,0)$-tensor $\nabla^g_\cdot\alpha_-$ is zero. Here  $\nabla^g$ denotes the Levi-Civita connection of $(M,g)$.
Equivalently, the Dirac current  $V_\phi$ of $\phi$ is a conformal Killing vector, i.e., we have $\mathcal{L}_{V_\phi}g=\lambda\cdot g$ for the Lie derivative
$\mathcal{L}$ and some function $\lambda\in C^\infty(M)$. In fact, this property of the  Dirac current of a conformal Killing spinor is well known
(cf.\ e.g.~\cite{B00})

On the other hand, we have also shown in \cite{LeiNCK} how to recover a  $\nabla^\cT$-parallel $q$-form tractor from its projecting component with respect to a choice of
a metric $g$ in the conformal class. Here,  in our case for the induced $2$-form tractor $\alpha_\Phi$ of a $\nabla^\cW$-parallel twistor
$\Phi=(\phi,\frac{\sqrt{2}}{n}\cD^g\phi)$, we have
\begin{gather}
\alpha_- = \alpha_\phi,\nonumber\\
\alpha_\mp=\frac{1}{n}d^*\alpha_\phi = \frac{2}{n}\mbox{Re}\langle\phi,\cD\phi\rangle_\cS,\nonumber\\
\alpha_+ = \Box\alpha_\phi =  \frac{-2}{n^2}\cdot \alpha_{\cD\phi},\label{eq3mal}
\end{gather}
where $d^*$ denotes the codif\/ferential and  $\Box$ is def\/ined as $\frac{-1}{n-2}(\Delta^g-{\rm tr}_g\Rho^g)$ with Bochner-Laplacian  $\Delta^g={\rm tr}_g(\nabla^g)^2$. The $2$-form
component of
$\alpha_\Phi$ is given with respect to  $X,Y\in TM$
by
\begin{equation}\label{eq1mal} \alpha_0(X,Y)  = \frac{1}{2}d\alpha_\phi(X,Y)  = \frac{1}{n}\mbox{Re}\langle X\wedge Y\cdot\phi,\cD\phi\rangle_\cS  .
\end{equation}

\section{A local geometric description}\label{Sec5}

In this section we derive a local geometric description of conformal Lorentzian manifolds admitting a conformal Killing spinor (outside of its  singular set).
There occur two possible geometries, the so-called {\it Brinkmann spaces} and the {\it static monopoles} with parallel spinors.
Our argument utilises tractor calculus in an essential way.

Let $(M^n,g)$, $n\geq 3$, be a connected, time-oriented Lorentzian spin manifold with non-trivial conformal Killing spinor   $\phi\in\Gamma(\cS)$ and  Dirac current
$V_\phi$.
We have ${\rm zero}(\phi)={\rm zero}(V_\phi)$, and~$V_\phi$ is outside of its zero set either null (i.e. $g(V_\phi,V_\phi)=0$) or timelike  (i.e.\ $g(V_\phi,V_\phi)<0$)
(cf.\ Lemma~\ref{lemorb}).
If $V_\phi$ is nowhere timelike on $M$, then we def\/ine the singular set of $\phi$ by \[{\rm sing}(\phi):={\rm zero}(\phi) .\] If $V_\phi$ is timelike in at least one point
of $M$, then we set \[{\rm sing}(\phi):={\rm zero}(g(V_\phi,V_\phi))  .\]
Recall that the corresponding $\nabla^\cW$-parallel twistor to $\phi$ is given with respect to $g$ by $\Phi=(\phi,\frac{\sqrt{2}}{n}\cD\phi)\in\Gamma(\cW)$ and the
induced
$2$-form
tractor
is $\alpha_\Phi=(\alpha_-,\alpha_0,\alpha_\mp,\alpha_+)$.

It is a matter of fact that  the $\SO(2,n)$-module $\frak{so}(2,n)\cong\Lambda^2({\RR^{2,n}}^*)$ decomposes under the adjoint action into the disjoint union of
$\SO(2,n)$-orbits.
A distinguished element of such an  $\SO(2,n)$-orbit can then be seen as a  normal form for all the  members of that orbit.
Note that any f\/iber $\cA_p$, $p\in M$, of the $2$-form tractor bundle $\cA$ on $(M,[g])$ is uniquely
identif\/ied with $\Lambda^2({\RR^{2,n}}^*)$ via the choice of a tractor frame  at $p\in M$ (i.e.\ an orthonormal basis of $\cT_p$).
Via such an identif\/ication any  $2$-form tractor in $\cA_p$ belongs to a uniquely determined  $\SO(2,n)$-orbit of $\Lambda^2({\RR^{2,n}}^*)$. If
two $2$-form tractors $A_p\in \cA_p$ and $A_q\in \cA_q$, $p,q\in M$, belong to the same orbit via an  identif\/ication with $\Lambda^2({\RR^{2,n}}^*)$,
we say  $A_p$ and  $A_q$ have the same orbit type.   In our situation, the $2$-form tractor
$\alpha_\Phi$, which belongs to the conformal Killing spinor $\phi$ on $(M,[g])$, is $\nabla^\cT$-parallel.
Hence, since the tractor connection $\nabla^\cT$ is induced by a Cartan connection
with values in  $\frak{so}(2,n)$, it is immediately clear that the $\SO(2,n)$-orbit type of $\alpha_\Phi$ at $p$ is constant all over  $M$.

Now let us assume that the zero set ${\rm zero}(\phi)$ of the non-trivial conformal Killing spinors $\phi$ is non-empty on $(M,g)$. Let $p\in {\rm zero}(\phi)$ be such
a zero of~$\phi$. It follows immediately from~(\ref{eq3mal}) and~(\ref{eq1mal}) that $d\alpha_\phi(p)=0$ and $d^*\alpha_\phi(p)=0$.
Also note that $\Box\alpha_\phi$
cannot vanish at $p$, since otherwise the  $2$-form tractor  $\alpha_\Phi$ had to be constant zero on~$M$. This shows that $\alpha_\Phi$ is given at $p\in M$ by
\begin{equation}\label{eq4}
\alpha_\Phi(p)  = s_+^\flat(p)\wedge\Box\alpha_\phi(p)  ,
\end{equation}
and this expression determines the orbit type of $\alpha_\Phi$ not only at $p$, but everywhere on $M$.
In fact, (\ref{eq4}) shows that the orbit type is represented by
a simple wedge product $l\wedge c$ of a null $1$-form $l\in{\RR^{2,n}}^*$ with another $1$-form $c\neq 0$.
Since this $1$-form $c$
corresponds to the Dirac current of the spinor $\cD\phi$, it has to be  causal. We conclude that the $2$-form tractor $\alpha_\Phi$,
which belongs via the twistor $\Phi$  to the conformal Killing spinor $\phi\in\Gamma(\cS)$ with zero at $p$, has two possible normal forms in  $\Lambda^2({\RR^{2,n}}^*)$
under the adjoint action of $\SO(2,n)$, namely either $l\wedge c$ with $c$ null or $l\wedge c$ with $c$ timelike!

Moreover, let us say that $V_\phi$ is {\it hypersurface orthogonal} if the $g$-orthogonal complement $V_\phi^\bot\subset TM$ of $V_\phi$ is an integrable distribution of
codimension
$1$ in $TM$ (over $M \backslash\, {\rm zero}(V_\phi)$).

\begin{lemma} \label{lem1} Let $\phi\in\Gamma(\cS)$ be a non-trivial conformal Killing spinor on a Lorentzian spin manifold $(M^n,g)$ of dimension $n\geq 3$ with
${\rm zero}(\phi)\neq
\varnothing$.
Then the
Dirac current $V_\phi$ of $\phi$ is hypersurface orthogonal on $M\backslash\, {\rm zero}(\phi)$.
\end{lemma}

\begin{proof} Since $\alpha_\Phi$ is simple, we have $\alpha_\Phi\wedge \alpha_\Phi{=}0$ on $(M,[g])$. In particular, this proves $\alpha_\phi\wedge d\alpha_\phi{=}0$
on  $M$. It follows directly by Frobenius' theorem that the annihilator $\{X\in TM|\, \iota_X\alpha_\phi=0\}$ is an integrable distribution of codimension $1$ on
$M\backslash\,
{\rm zero}(\phi)$. The annihilator is the $g$-orthogonal complement of $V_\phi$ in $TM$.
\end{proof}

It is a well known fact that any hypersurface orthogonal, conformal Killing vector f\/ield $V$, which does not change its causal type, is (at least locally) parallel with
respect to
the Levi-Civita connection $\nabla^{\tilde{g}}$ of some metric $\tilde{g}$ in the given conformal class $c=[g]$ of a Lorentzian manifold $(M,g)$
(cf.\ e.g.~\cite{L07b}).  We call a Lorentzian metric $\tilde{g}$, which admits a parallel null vector f\/ield, a~{\it Brinkmann metric}.

\begin{proposition} \label{prop1} Let $\phi\not\equiv 0$ be a conformal Killing spinor on a connected Lorentzian spin manifold $(M^n,g)$ of dimension $n\geq 3$  with
${\rm zero}(\phi)\neq
\varnothing$. Then the set $M\backslash\, {\rm sing}(\phi)$ is dense in $M$, and
for any $q\not\in {\rm sing}(\phi)$ there exists a neighbourhood $U_q\subset M\backslash\, {\rm sing}(\phi)$ and a function $f$ on $U_q$ such that
the rescaled spinor $\tilde{\phi}:=e^{-f/2}\phi$ and the  Dirac current $V_\phi=V_{\tilde{\phi}}$ are parallel with respect
to the metric $\tilde{g}:=e^{2f}g$ on $U_q$. Two cases are possible.
\begin{enumerate}\itemsep=0pt
\item[$1.$] The  Dirac current $V_\phi$ is timelike in some point of $M$ and
the Lorentzian metric $\tilde{g}$ on $U_q$  is isometric to a static monopole
$-dt^2 + h$,
where $V_\phi=\partial t$ and $h$ is a Ricci-flat Riemannian metric with parallel spinor.
\item[$2.$] The  Dirac current $V_\phi$ is null on  $M\backslash\, {\rm zero}(\phi)$ and $\tilde{g}$ on $U_q$  is a Brinkmann metric with parallel spinor.
\end{enumerate}
\end{proposition}

\begin{proof}
First, let us assume that $f$ is some function on an open subset $U$ in $M$ such that  $V_\phi$ is parallel with respect to
the Levi-Civita connection $\nabla^{\tilde{g}}$ of $\tilde{g}=e^{2f}g$. Then the   $2$-form tractor~$\alpha_\Phi$ to $\tilde{\phi}=e^{-f/2}g$
must be given on $U$ with respect to $\tilde{g}$ by $s_-^\flat\wedge\alpha_{\tilde{\phi}}+s_+^\flat\wedge\Box^{\tilde{g}}\alpha_{\tilde{\phi}}$.  However,
since the orbit type of  $\alpha_\Phi$ is $l\wedge c$ with $c$ causal, we can conclude $\Box^{\tilde{g}}\alpha_{\tilde{\phi}}=0$
and $\alpha_\Phi=s_-^\flat\wedge\alpha_{\tilde{\phi}}$. This shows $\mathcal{D}^{\tilde{g}}\tilde{\phi}=0$ on $U$. Hence  $\tilde{\phi}$ is a parallel spinor
on $(U,\tilde{g})$.

Next we see that if $V_\phi$ is timelike at $p\in M$, then $V_\phi$
cannot be null on an open subset $W$ of $M$. (Otherwise the $1$-form $c$ of the orbit type of  $\alpha_\Phi$ would be null on $W$ and timelike at $p$.)
Hence, since $M\backslash\, {\rm zero}(\phi)$ is dense in $M$, we can also conclude that $M\backslash\, {\rm sing}(\phi)$ is dense in $M$.

It remains to describe  the local conformal geometry of $g$ of\/f the singularity set on $M$. First, let us assume that
$V_\phi$ is timelike on  $M\backslash\,  {\rm sing}(\phi)$. In this case we can rescale to $\tilde{g}:=(g(V_\phi,V_\phi))^{-2}\cdot g$ on  $M\backslash\,  {\rm sing}(\phi)$.
Since $V_\phi$ is parallel with respect to $\tilde{g}$, there exits locally a coordinate $t$ such that $\tilde{g}$ is given by $-dt^2+h$ with some Riemannian
metric $h$. It is well known in this case that the  $\tilde{g}$-parallel spinor $\tilde{\phi}$ induces a parallel spinor for the Riemannian metric $h$. In
particular, $h$ is Ricci-f\/lat and  $\tilde{g}=-dt^2+h$ is a static monopole metric.

In the other case, the Dirac current $V_\phi$ is null and locally parallel with respect to some metric $\tilde{g}=e^{2f}g$ on an appropriate neighbourhood $U_q$ of any
point
$q\in M\backslash\, {\rm zero}(\phi)$. This means that $\tilde{g}$ is a Brinkmann metric. The parallel spinor is  $\tilde{\phi}:=e^{-f/2}\phi$ on $U_q$.
\end{proof}

\begin{remark} We want to emphasis here that the  arguments in this section essentially rely on tractor calculus. The assumption ${\rm zero}(\phi)\neq \varnothing$
in  Lemma \ref{lem1} and
Proposition \ref{prop1} determines via tractor calculus the orbit type of $\alpha_\Phi$, and hence the conformal geometry not only next to
the zero set, but everywhere on $M\backslash\,  {\rm sing}(\phi)$. Note that, in general, a conformal Killing spinor without zero on a Lorentzian spin manifold
is neither conformally related to a parallel spinor nor is its
Dirac current hypersurface orthogonal.  In fact, certain  conformal Killing spinors on Lorentzian Einstein--Sasaki spaces do not have this property
(cf.\ e.g.~\cite{B00}).

Also note that an analogous result in Riemannian geometry is well known. Any conformal Killing spinor $\phi$ with zero on  a Riemannian manifold $(N,h)$
is conformally related to a parallel spinor of\/f its zero set. A proof of this result does not need tractor calculus. One can rescale
the Riemannian metric $h$ with $\|\phi\|^{-4}$, where $\|\phi\|$ is the spinor norm, and argue that the rescaled spinor $\tilde{\phi}=\frac{1}{\|\phi\|}\phi$ is parallel
(cf.~\cite{BFGK91}).
In Lorentzian
geometry this kind of proof
does not
work, alone for the simple reason that
the Hermitian product $\langle\cdot,\cdot\rangle_\cS$ on the spinor bundle is indef\/inite. It seems that a proof with tractor methods is
unavoidable in the Lorentzian case.
\end{remark}

\section{The shape of the zero set}\label{Sec6}

Finally, we discuss here the shape of the zero set  ${\rm zero}(\phi)$ and the singular set ${\rm sing}(\phi)$
of a~conformal Killing spinor $\phi$ on a Lorentzian spin manifold. The results of this section are based on  arguments from~\cite{L99} and~\cite{L01}.

Let $(M^n,g)$, $n\geq 3$, be a Lorentzian spin manifold with spinor bundle $\cS$. The spinor squaring has the following useful properties.

\begin{lemma}[\cite{L01,BL04}] \label{lemorb} Let $\phi\in \cS_p$ be a non-trivial spinor at $p\in M$ with spinor square $V_\phi\in T_pM$. Then
\begin{enumerate}\itemsep=0pt
\item[$1.$] $V_\phi$ is causal.
\item[$2.$] $X\cdot\phi=0$ if and only if $V_\phi$ is null and $X\in\RR V_\phi$.
\end{enumerate}
\end{lemma}

In the following, we denote by $\gamma_p:I=(a,b)\to M$, $a<0<b\in\RR$,  a geodesic of the Lorentzian manifold $(M,g)$ with $\gamma_p(0)=p$.
We call a geodesic  $\gamma_p$ {\it maximal} on $M$ if the interval~$I$ is the maximal domain of def\/inition.
For an  arbitrary point $p$ of the  Lorentzian manifold $(M,g)$, we call the set $L_p$ of those points in $M$, which can be joined with~$p$ by a
smooth null
geodesic,
the {\it geodesic null cone} with origin~$p$.
We have the following general result from~\cite{L99} (see also~\cite{L01,F07}).

\begin{lemma} \label{lemalt} Let $V$ be a conformal vector field on a Lorentzian manifold $(M^n,g)$ with  $\nabla^gV(p)=0$ for all  $p\in {\rm zero}(V)$. Then
\begin{enumerate}\itemsep=0pt
\item[$1.$] There exists for any $p\in {\rm zero}(V)$ an open  neighbourhood $U_p\subset M$ such that  ${\rm zero}(V)\cap U_p$ is contained in the image ${\rm Im}(\gamma_p)$ of a smooth null
geodesic $\gamma_p$.
\item[$2.$] For any point $q\in L_p$ of the null cone of $p\in {\rm zero}(V)$, the vector $V(q)\in T_qM$ is a multiple of the tangent vector to a null geodesic, which
joins $p$ and $q$. In particular, any $V(q)\neq 0$ is null for $q\in L_p$ and $p\in {\rm zero}(V)$.
\end{enumerate}
\end{lemma}

Now let $\phi\in\Gamma(\cS)$ be a conformal Killing spinor on $(M^n,g)$ with ${\rm zero}(\phi)\neq\varnothing$. Note that the relations (\ref{eq3mal}) and (\ref{eq1mal}) show
$(\nabla^gV_\phi)(p)=0$ for any $p\in {\rm zero}(\phi)$, i.e., Lemma~\ref{lemalt} is applicable for the Dirac current $V_\phi$ of $\phi$.

\begin{lemma}  \label{lemmehr} Let $\phi\in\Gamma(\cS)$ be a conformal Killing spinor with zero on a Lorentzian spin manifold $(M^n,g)$,
and let $\gamma_p$ be a geodesic on $M$ with $\gamma_p(0)=p$.
\begin{enumerate}\itemsep=0pt
\item[$1.$] If $\gamma_p'(0)\cdot \cD\phi(p)=0$ then ${\rm Im}(\gamma_p)\subset {\rm zero}(\phi)$.
\item[$2.$] If  $\gamma_p'(0)\cdot \cD\phi(p)\neq 0$ then there exists a neighbourhood $U_p$ of $p$ in $M$ such that $p$ is the only zero of $\phi$ in the image
of $\gamma_p$  in $U_p$.
\end{enumerate}
\end{lemma}

\begin{proof} Let $\mathcal{B}(t)=(b_1(t),\ldots,b_m(t))$ be a $\nabla^\cS$-parallel complex spinor frame along  $\gamma_p(t)$. (The complex rank of
$\cS$ is $m:=2^{[\frac{n}{2}]}$.) We set $u_i(t):=\langle\phi(\gamma_p(t)),b_i(t)\rangle_\cS$ for $i=1,\ldots,m$. Then we have
\begin{gather*}
u_i' = -\frac{1}{n}\langle\gamma_p'\cdot \cD\phi,b_i\rangle_\cS\qquad\mbox{and}\qquad
u_i'' = -\frac{1}{2}\langle\gamma_p'\Rho(\gamma_p')\cdot \phi,b_i\rangle_\cS  ,
\end{gather*}
i.e., the complex vector $U(t):=\left(\begin{array}{c}u_1(t)\\\vdots\\ u_m(t)\end{array}\right)$ satisf\/ies the ordinary  dif\/ferential equation $U''=-\frac{1}{2}C\cdot U$,
where $C\in M(m,\CC)$ is the complex matrix of the endomorphism $\phi\in\cS\mapsto \gamma_p'\Rho(\gamma_p')\cdot\phi\in\cS$ with respect to  $\mathcal{B}(t)$.
By assumption, we  have $U(0)=0$. Moreover, $U'(0)=0$ if and only if $\gamma_p'(0)\cdot\cD\phi(p)=0$. Obviously, the general  theory of ordinary dif\/ferential equations
implies the two statements of Lemma~\ref{lemmehr}.
\end{proof}

We obtain the following result about the shapes of ${\rm zero}(\phi)$ and ${\rm sing}(\phi)$.

\begin{proposition} \label{propfast}  Let $\phi\in\Gamma(\cS)$ be a non-trivial conformal Killing spinor with ${\rm zero}(\phi)\neq \varnothing$ on a~connected Lorentzian spin
manifold
$(M^n,g)$.
\begin{enumerate}\itemsep=0pt
\item[$1.$] If $V_\phi$ is null on $M\backslash\, {\rm zero}(\phi)$, then   ${\rm sing}(\phi)={\rm zero}(\phi)$ is a countable union of isolated images of smooth  maximal null geodesics.
\item[$2.$]  If $V_\phi$ is timelike at some point $p\in M$, then ${\rm zero}(\phi)$ is a countable union of isolated points.
The singular set ${\rm sing}(\phi)$ contains
all geodesic null cones $L_p$ with origin $p\in {\rm zero}(\phi)$. In some neighbourhood $U_p\subset M$ of any $p\in {\rm zero}(\phi)$, the singular set ${\rm sing}(\phi)$ equals
the null cone $U_p\cap
L_p$.
\end{enumerate}
\end{proposition}

\begin{proof} (1) First, let us assume that $V_\phi$ is null on $M\backslash\, {\rm zero}(\phi)$. Then we know from (\ref{eq4}) that $V_{\cD\phi}(p)$ is null in $p\in
{\rm zero}(\phi)$.
By Lemma \ref{lemorb} and Lemma \ref{lemmehr} (1),  we know that the image of the maximal null geodesic $\gamma_p$ with  $\gamma_p'(0)=V_{\cD\phi}(p)$ is contained in
${\rm zero}(\phi)$. Now,
Lemma
\ref{lemalt} implies that in a neighbourhood $U_p$ of $p$ the zero set  ${\rm zero}(\phi)$ is equal to $U_p\cap {\rm Im}(\gamma_p)$.
However, since this is true for any zero $p$ of $\phi$,
the assertion (1) of Proposition \ref{propfast} follows on a manifold with countable basis.

(2)  Let $V_\phi$ be timelike at some point of $M$. Then from (\ref{eq4}) we know that $V_{\cD\phi}(p)$ is timelike at $p\in
{\rm zero}(\phi)$. Hence, by Lemma \ref{lemorb},  $\gamma_p'(0)\cdot \cD\phi(p)=0$  is
impossible
for
any
geodesic through $p\in
{\rm zero}(\phi)$, and
Lemma \ref{lemalt} implies that the zero $p$ is an isolated point. This is true for any $p\in {\rm zero}(\phi)$.

It remains to prove the statements about the singular set
${\rm sing}(\phi)$. Obviously, Lemma   \ref{lemalt} (2) implies that $L_p$ is contained in ${\rm sing}(\phi)$ for any $p\in {\rm zero}(\phi)$. Now, let us assume that $V_\phi(q)$ is
null at a point $q$ in a small
neighbourhood $U_p$ of $p$ with $q\not\in L_p$. Since $V_\phi$ is a conformal vector f\/ield, the integral curve to $V_\phi$ through $q$ is a null curve in $M$.
In fact, the integral curve to~$V_\phi$ through~$q$ is a null (pre-)geodesic (cf.~\cite[p. 63]{L01}). However, a straightforward argument  as in~\cite{L99} shows that
any
maximal
null
geodesic, which runs through a point $q$ that  is {\it very close} to~$p$, has to intersect the null cone $L_p$ of $p$ in a point $\ell$. This intersection is transversal
to
$L_p$.
Hence the Dirac current $V_\phi$  has to be zero at $\ell$ (since $V_\phi$ is tangent to two transverse null geodesics running through $\ell$). However, if we
consider a
point
$q$, which is  arbitrary  close to~$p$, then the
intersection
point $\ell$ is arbitrary close to $p$ as well. But $p$ is an isolated zero of $\phi$. This implies that there cannot exist  a point $q$ arbitrary close to $p\in
{\rm zero}(\phi)$ of\/f the null cone $L_p$,  where $V_\phi$ is null.
\end{proof}

In summary with Proposition \ref{prop1} we obtain our main result.

\begin{theorem} \label{The}  Let $\phi\in\Gamma(\cS)$ be a non-trivial conformal Killing spinor with ${\rm zero}(\phi)\neq \varnothing$ on a~connected Lorentzian spin
manifold $(M^n,g)$ of dimension $n\geq 3$. Then ${\rm zero}(\phi)$ consists either~of
\begin{enumerate}\itemsep=0pt
\item[$1.$] isolated images of null geodesics and off the zero set the metric $g$ is locally conformally equivalent to a Brinkmann metric with parallel spinor,  or
\item[$2.$]   isolated points and  off the singular set ${\rm sing}(\phi)$  the metric $g$ is conformally equivalent to a~static monopole $-dt^2+h$, where
$h$ is a Riemannian metric with parallel spinor.
\end{enumerate}
\end{theorem}

\begin{remark} Both cases of Theorem \ref{The} do occur. In fact, any conformal Killing spinor $\phi$ on the Minkowski space $\RR^{1,n-1}$ with a zero at the origin
is explicitly given by \[\phi:\ p\in \RR^{1,n-1}  \mapsto   p\cdot S\qquad \mbox{with}\quad S\in\mathbb{\Delta}_{1,n-1}  .\]   For those
$S\in\mathbb{\Delta}_{1,n-1}$,
whose spinor square is null, the zero set of $\phi$ is a null line in $\RR^{1,n-1}$ (in direction of the spinor square).  For those
$S\in\mathbb{\Delta}_{1,n-1}$,
whose spinor square is timelike,  the zero set consists only of the origin, and the null cone of  $\RR^{1,n-1}$ is the singular set.

Moreover, note that in~\cite{L07a} we have explicitly constructed a conformal Killing spinor  with an isolated zero on a Lorentzian spin manifold, which is not
conformally f\/lat in any neighbourhood of the  isolated zero. However, this Lorentzian metric is not smooth, but only $C^1$-dif\/ferentiable, along the singular set.
On the other hand, in \cite{F07} it is stated that any analytic
Lorentzian metric, admitting a causal conformal vector f\/ield with a zero, is
globally conformally f\/lat.
\end{remark}

\pdfbookmark[1]{References}{ref}
\LastPageEnding

\end{document}